\documentclass[11pt]{article}
\usepackage[T1]{fontenc}
\usepackage{lmodern,amsmath,amssymb,amsthm,mathtools}
\usepackage{enumitem}
\usepackage[margin=28mm]{geometry}
\usepackage{microtype}
\usepackage[colorlinks=true,linkcolor=blue,citecolor=blue,urlcolor=blue]{hyperref}

\newtheorem{theorem}{Theorem}
\newtheorem{lemma}[theorem]{Lemma}
\newtheorem{conjecture}[theorem]{Conjecture}
\newtheorem{question}[theorem]{Question}

\newcommand{\vp}{\varphi}
\newcommand{\Z}{\mathbb Z}

\newcommand{\ora}{\overrightarrow}
\DeclareMathOperator{\dist}{dist}

\title{The minimum number of cycles of length $0$ modulo $k$ in graphs with no proper $k$-coloring}

\author{
Ilkyoo Choi\thanks{Department of Mathematics, Hankuk University of Foreign Studies, Yongin-si, Gyeonggi-do, Republic of Korea.
 E-mail: \texttt {ilkyoo@hufs.ac.kr},
  Discrete Mathematics Group, Institute for Basic Science (IBS), Daejeon, Republic of Korea, and Center for AI and Natural Sciences, Korea Institute for Advanced Study (KIAS), Seoul, Republic of Korea.
Supported by the Hankuk University of Foreign Studies Research Fund, the National Research Foundation of Korea (NRF) grant funded by the Korea government (MSIT) (RS-2025-23324220), Institute for Basic Science (IBS-R029-C1), and the Korea Institute for Advanced Study (KIAS) grant funded by the Korea government.
}
\and
Kamil Sieniawski\thanks{Independent researcher, Pozna\'n, Poland.  E-mail: \texttt{research@sieniawski.dev}}
}

\begin{document}
\maketitle

\begin{abstract}
Consider an integer $k\geq 3$.
Toft and Tuza observed that every graph with no proper $k$-coloring has a cycle of length $0$ modulo $k$.
Brewster, McGuinness, Moore, and Noel conjectured the following quantitative strengthening: every graph with no proper $k$-coloring has at least $\frac{(k+1)!}{2k}$  cycles whose lengths are divisible by $k$. 
The conjecture is sharp,  as demonstrated by the complete graph on $k+1$ vertices.
Kim and Picollelli confirmed the conjecture for $k=3$.
We establish the conjecture for every $k\geq 3$. 
In other words, we determine the minimum number of cycles of length
$0$ modulo $k$ in a graph with no proper $k$-coloring. 
\end{abstract}

\section{Introduction}

All graphs considered in this paper are finite and simple, which means no loops and no parallel edges. 
Let $K_n$ denote the complete graph on $n$ vertices, and an \emph{$\ell$-cycle} is a cycle of length $\ell$.
Given a graph $G$, a \emph{proper $k$-coloring} of  $G$ is a function $\vp:V(G)\to\{0,\ldots,k-1\}$ such that $\vp(x)\ne\vp(y)$ for every edge $xy$. 
The minimum $k$ for which $G$ has a proper $k$-coloring is the \emph{chromatic number} of $G$, denoted $\chi(G)$. 

A trivial observation is that every graph with no proper $2$-coloring  has an odd cycle, which is a cycle of length $1$ modulo $2$.
Tuza~\cite{Tuza1992} extended this to every integer $k\geq 3$: every graph with no proper $k$-coloring has a cycle of length $1$ modulo $k$.
In \cite{Tuza1992}, it is also stated that Toft and Tuza observed that for $k\geq 3$, every graph with no proper $k$-coloring  has a cycle of length $0$ modulo $k$.
Note that the bound on $k$ is sharp as an odd cycle has neither a proper $2$-coloring nor a cycle of length $0$ modulo $2$.
More generally, Chen, Ma, and Zang~\cite{ChenMaZang2015} proved a result implying that for $k\ge3$ and $r\in\{0,\ldots,k-1\}\setminus\{2\}$, every graph with no proper $k$-coloring has a cycle of length $r$ modulo $k$.
The exception of $r=2$ is necessary as $K_{k+1}$ has neither a proper $k$-coloring nor a  cycle of length $2$ modulo $k$.
For $r=2$, Diwan, Kenkre, and Vishwanathan~\cite{2013DiKeVi} showed that every graph with no proper $(k+1)$-coloring has a cycle of length $2$ modulo $k$. 
Resolving a conjecture by Sudakov and Verstra\"ete~\cite{2017SuVe}, Gao, Huo, Liu, and Ma~\cite{2022GaHuLiMa} proved that every graph with no proper $k$-coloring has $k-1$ cycles of consecutive lengths. 
See also~\cite{2021GaHuMa}.

A graph with no proper $k$-coloring has a subgraph with minimum degree at least $k$, otherwise, the graph is $(k-1)$-degenerate and thus has a proper $k$-coloring. 
Dean~\cite{1991Dean} conjectured that for $k\geq 3$, every graph  with minimum degree at least $k$ has a cycle of length $0$ modulo $k$. 
This conjecture was proved  for $k=3$ by Chen and Saito~\cite{ChenSaito1994} and for $k=4$ by Dean, Lesniak, and Saito~\cite{DeanLesniakSaito1993}.
Recently, Luo, Ma, and Zhao~\cite{LuoMaZhao2026} proved Dean's conjecture for every $k\ge6$, leaving only the case $k=5$ open. 

Once existence is guaranteed, a natural followup question is a quantitative question: given a graph with no proper $k$-coloring, what is the guaranteed number of cycles of length $0$ modulo $k$? 
For $k\geq 3$, the complete graph $K_{k+1}$ has no proper $k$-coloring and has $(k+1)\frac{(k-1)!}{2}
=\frac{(k+1)!}{2k}$ cycles of length $0$ modulo $k$. 
Note that the only possible cycles of length $0$ modulo $k$ in $K_{k+1}$ are $k$-cycles.
Brewster, McGuinness, Moore, and Noel~\cite{BMMN2016} conjectured
that among all graphs with no proper $k$-coloring, $K_{k+1}$ minimizes the number of cycles of length $0$ modulo $k$.

\begin{conjecture}[Brewster, McGuinness, Moore, and Noel~\cite{BMMN2016}]
\label{conj:BMMN}
For every $k\geq 3$, every graph with no proper $k$-coloring has at least $\frac{(k+1)!}{2k}$ cycles of length $0$ modulo $k$. 
\end{conjecture}

Brewster, McGuinness, Moore, and Noel~\cite{BMMN2016} proved the following factorial lower bound,
which is smaller than  the conjectured bound by a factor of $k+1$.

\begin{theorem}[Brewster, McGuinness, Moore, and Noel~\cite{BMMN2016}]
\label{thm:BMMN}
For every $k\geq 3$, every graph with no proper $k$-coloring has at least $\frac{(k-1)!}{2}$ cycles  of length $0$ modulo $k$. 
\end{theorem}

Strengthening results in~\cite{Tuza1992}, Moore and West~\cite{MooreWest2021} proved further counting results regarding the number of cycles of certain lengths modulo $k$ through a given edge. 
They also proved that for $k\geq3$, every graph with no proper $k$-coloring has  $\frac{k!}2$ cycles of length $1$ modulo $k$.
Note that this bound is sharp as again demonstrated by $K_{k+1}$. 
For the number of cycles of length $0$ modulo $k$ in a graph with no  proper $k$-coloring, their results guarantee $\frac{(k-1)!}2$ cycles, recovering the result in Theorem~\ref{thm:BMMN}.

Kim and Picollelli~\cite{KimPicollelli2024} confirmed Conjecture~\ref{conj:BMMN} for $k=3$. 
We establish Conjecture~\ref{conj:BMMN} for every $k\geq 3$. 

\begin{theorem}\label{thm:main}
For every $k\ge3$, every graph with no proper $k$-coloring has at least $\frac{(k+1)!}{2k}$ cycles of length $0$ modulo $k$. 
\end{theorem}

Recall that the bound of $\frac{(k+1)!}{2k}$ is tight by $K_{k+1}$.

Our proof uses the approach of Brewster, McGuinness, Moore, and Noel~\cite{BMMN2016}.
Relative to a vertex ordering, an $ak$-cycle is \emph{good} if it can be traversed with at most $a$ backward edges, where an edge is backward when traversed from a later vertex to an earlier vertex in the vertex ordering.
In Section~\ref{sec:cycles}, we show that if $H$ has no proper $k$-coloring but $H-v$ has a proper $k$-coloring, then every ordering induced by a proper $k$-coloring of $H-v$, with $v$ placed last, has at least $k+1$  good cycles.
In Section~\ref{sec:count}, we show that a fixed cycle is good for at most $2k$ permutations of the colors of a fixed proper $k$-coloring.
Finally, we in prove Theorem~\ref{thm:main} by combining these bounds in Section~\ref{sec:conclusion}.
We also suggest some future directions. 

 \section{Counting good cycles}
 \label{sec:cycles}

In this section, we show that certain vertex orderings of a graph with chromatic number at least $k+1$ yield at least $k+1$ good cycles.
We begin with a recoloring lemma.

Fix $k\ge3$. 
Let $H$ be a graph. 
Given a labeling $\vp:V(H)\to\Z$, an edge $xy$ is \emph{$\vp$-monochromatic} if $\vp(x)\equiv\vp(y)\pmod k$. 
Let $D(H, \vp)$ be the digraph on $V(H)$ where $\ora{xy}$ is an arc if and only if $xy\in E(H)$ and $\vp(y)-\vp(x)\equiv1\pmod k$.
Thus a directed path of length $r$ increases the label by $r$, modulo $k$.
Clearly, $D(H, \vp)$ has no loops and no directed $2$-cycles.

The following lemma is a version of the recoloring argument used by Moore and West~\cite{MooreWest2021};
see also Kim and Picollelli~\cite{KimPicollelli2024}.

\begin{lemma}\label{lem:recolor}
For $k\geq3$, let $H$ be a graph with $\chi(H)>k$, and let $\vp:V(H)\to\Z$ be a labeling of $H$. 
For each $\vp$-monochromatic edge $e$, let the two endpoints be $t(e)$
and $h(e)$.
Then there exist $\vp$-monochromatic edges $e_1, e_2$ such that $t(e_1)$ can reach $h(e_2)$ in $D(H,\vp)$.
\end{lemma}
\begin{proof}
Let $M$ be the set of $\vp$-monochromatic edges of $H$.
Since $H$ has no proper $k$-coloring, $M$ is nonempty;
otherwise, replacing each label $\vp(v)$ with $\vp(v)$ modulo $k$ is a proper $k$-coloring of $H$, a contradiction. 
Let $R$ be the set of vertices reachable in $D(H,\vp)$ from $t(e)$ for some $e\in M$. 
Note that $t(e)\in R$ for every $e\in M$.

Suppose to the contrary that $h(e)\notin R$ for every $e\in M$.
Define a labeling $\psi:V(H)\to\{0,\ldots,k-1\}$ where $\psi(u)\equiv
\begin{cases}
\vp(u)+1 & \pmod k, u\in R,\\
\vp(u) & \pmod k, u\notin R.
\end{cases}$
Every edge $e\in M$ is properly colored by $\psi$, since $t(e)\in R$ and $h(e)\notin R$.
Assume $xy$ is monochromatic by $\psi$.
Thus $xy\notin M$ and without loss of generality let $x\in R$ and $y\notin R$.
Since $\psi(x)=\psi(y)$, we know $\vp(y)-\vp(x)\equiv1\pmod k$.
Thus, $\ora{xy}$ is an arc of $D(H,\vp)$.
Since $x\in R$, this implies $y\in R$, a contradiction. 
Thus $\psi$ is a proper $k$-coloring of $H$, a contradiction.
\end{proof}

Let $H_\prec$ denote the graph $H$ with an ordering $\prec$ on its vertices. 
A traversal of an edge $xy\in E(H_\prec)$ from $x$ to $y$ is \emph{forward} if $x\prec y$ and \emph{backward} otherwise.
For a cycle $C$ in $H_\prec$, let $\vec C$ denote $C$ with one of the two cyclic traversals fixed. 
We call $\vec C$ \emph{good} if $C$ has length $ak$ and $\vec C$ has at most $a$ backward edges.
We also call $C$ \emph{good} if one of its two cyclic traversals is good. 
Note that the two cyclic traversals cannot be both good.
One has $b$ backward edges and the other has $ak-b$ backward edges. 
If both were good, then $ak\leq a+b\leq 2a$, which contradicts $k\ge3$.

\begin{lemma}\label{lem:orders}
For $k\ge3$, let $H$ be a graph with $\chi(H)>k$. 
Suppose $H$ has a vertex $v$ such that $H-v$ has a proper $k$-coloring $\vp:V(H-v)\to\{0, \ldots, k-1\}$. 
Let $\prec$ be a vertex ordering of $H$ obtained by ordering the vertices of $H-v$ in nondecreasing order of their colors under $\vp$, and placing $v$ last.
Then $H_\prec$ has at least $k+1$ distinct good cycles.
\end{lemma}
\begin{proof}
Extend $\vp$ to $V(H)$ by setting $\vp(v)=k-1$.
Let $Z=\{z\in N_H(v):\vp(z)=k-1\}$.
Since $H$ has no proper $k$-coloring, $Z$ is nonempty.
The only $\vp$-monochromatic edges of $H$ are $zv$ for $z\in Z$.
Let $D=D(H, \vp)$.

We prove the lemma in four steps.
First, we define $k+1$ new labelings $\vp_i$ obtained by modifying $\vp$.
Second, to the endpoints of each $\vp_i$-monochromatic edge, we assign $t(e)$ and $h(e)$.
Third, using Lemma~\ref{lem:recolor}, we find a good cycle $C_i$ for each $i\in\{1,\ldots,k+1\}$ containing a specific vertex $t_i$. 
Finally, we show that these $k+1$ good cycles are pairwise distinct by these $t_i$'s.

\begin{enumerate}[
    wide=\dimexpr\parindent/2\relax,
    label=\emph{{Step \arabic*.}},
]

\item \emph{New labelings $\vp_i$.}

For a vertex $u$, let $\dist(u)$ denote the length of a shortest directed path from $v$ to $u$ in $D$, and let $\dist(u)=\infty$
if $u$ is not reachable from $v$.
Note that each arc increases the label by $1$, modulo $k$.
For $z\in Z$, every finite $\dist(z)$ is a positive multiple of $k$ since $\vp(z)=\vp(v)$.
In particular, $\dist(z)\ge k$.

For $i\in\{0,\ldots,k+1\}$, let
$S_i=\{u\in V(H):\dist(u)\leq i-1\}$.
Note that $S_0=\emptyset$ and $S_1=\{v\}$. 
Let $T_i$ be the set of vertices that can reach some vertex of $Z$ in $D-S_i$.
We remark that paths of length $0$ are allowed.
Since $S_0\subseteq S_1\subseteq\cdots\subseteq S_{k+1}$, we know $T_0\supseteq T_1\supseteq\cdots\supseteq T_{k+1}$.
Each vertex $u\in T_i$ satisfies $\dist(u)\ge i$ and has a directed path to $Z$ entirely inside $T_i$.
For $z\in Z$, we have $z\in T_i$ if and only if
$\dist(z)\ge i$.
Since $\dist(v)=0$, we have $v\notin T_1$, in particular, $v\notin T_i$ for $i\geq 1$. 

For $i\in\{1,\ldots,k+1\}$, define a new labeling $\vp_i:V(H)\to\Z$ by decreasing the label of each vertex in $T_i$ by exactly $1$.
Note that vertices in $T_0$ receive a label of $-1$ under $\vp_i$. 
Let $D_i=D(H, \vp_i)$. 

\item \emph{Assigning $h(e)$ and $t(e)$ to endpoints of $\vp_i$-monochromatic edges.}

For each $\vp_i$-monochromatic edge $e$ of $H$, we assign $h(e)$ and $t(e)$ to its endpoints as follows. 

First, let $e=xy$ be a newly $\vp_i$-monochromatic edge.
Up to renaming $x$ and $y$ if necessary, $\ora{xy}$ is an arc of $D$ where $x\notin T_i$, $y\in T_i$, and $\vp(y)-\vp(x)\equiv1\pmod k$. 
Let $(t(e), h(e))=(x,y)$. 

We claim $\dist(x)+1=i=\dist(y)$ and $x\in T_{i-1}\setminus T_i$.
Indeed, let $P$ be a directed path from $y$ to $Z$ entirely in $T_i$, which exists since $y\in T_i$. 
If $\dist(x)\ge i$, then $\ora{xy}$ followed by $P$ shows that $x\in T_i$, a contradiction, so $\dist(x)\leq i-1$.
Since $y\in T_i$, we know $i\le\dist(y)\le\dist(x)+1\le i$, so $\dist(x)=i-1$ and $\dist(y)=i$.
Moreover, $\ora{xy}$ followed by $P$ is a directed path from $x$ to $Z$ in $D-S_{i-1}$, so $x\in T_{i-1}$.

Now, consider an originally $\vp$-monochromatic edge $zv$ of $H$, so $z\in Z$.
If $z\in T_i$, then $\vp_i(v)-\vp_i(z)=1$, so $\ora{zv}$ is an arc of $D_i$. 
Otherwise, $z\notin T_i$, so $\ora{zv}$ remains $\vp_i$-monochromatic and $\dist(z)<i\le k+1$.
Let $t(zv)=z$ and $h(zv)=v$. 
Since $\dist(z)$ is a positive multiple of $k$, we know $\dist(z)=k$ and $i=k+1$, so $z\in T_k\setminus T_{k+1}$.

Consequently, every $\vp_i$-monochromatic edge $e$ has
\begin{equation}\label{eq:tails}
\dist(t(e))=i-1,\qquad
t(e)\in T_{i-1}\setminus T_i,
\end{equation}
and either $h(e)\in T_i$ or $i=k+1$ anad $h(e)=v$.

\item \emph{Good cycles for $\vp_i$.}

We first show that every directed cycle in $D_i$ is good for  $\prec$.
Note that $i\geq 1$. 

Let $xy$ be an edge of $H_\prec$ with  $x\prec y$.
If $xy$ is not $\vp$-monochromatic, then $1\le\vp(y)-\vp(x)\le k-1$.
Since the label of an endpoint may  change by exactly $1$ under  $\vp_i$, we know $0\le\vp_i(y)-\vp_i(x)\le k$.
Otherwise, $xy$ is $\vp$-monochromatic, so $y=v$, $\vp(x)=k-1$, and $x\prec y$.
Since $v\notin T_i$, so know $\vp_i(y)-\vp_i(x)\in\{0,1\}$.
Hence, $0\le\vp_i(y)-\vp_i(x)\le k$ for every edge $xy$ with $x\prec y$.
Since every arc $\ora{xy}$ of $D_i$ satisfies $\vp_i(y)-\vp_i(x)\equiv1\pmod{k}$, we know $\vp_i(y)-\vp_i(x)=1$ if $x\prec y$ and $\vp_i(y)-\vp_i(x)=1-k$ if $y\prec x$. 

Now let $\vec C$ be a directed $\ell$-cycle in $D_i$ with $b$ backward edges.
Since $0=\sum_{\ora{xy}\in A(\vec C)}(\vp_i(y)-\vp_i(x))=(\ell-b)+b(1-k)=\ell-kb$, we know $\ell=kb$, so $\vec C$ is good.

We next show that $D_i$ contains a directed cycle through some vertex of $T_{i-1}\setminus T_i$.
Recall the assignment of $h(e)$ and $t(e)$ for each $\vp_i$-monochromatic edges from \emph{Step~2}. 
By Lemma~\ref{lem:recolor}, there are $\vp_i$-monochromatic edges $e_1, e_2$ such that $t(e_1)$ can reach $h(e_2)$ in $D_i$.
Let $(t_i, h_i)=(t(e_1), h(e_2))$.

We claim that $h_i$ can reach $t_i$ in $D_i$ through $v$.
We first show $h_i$ can reach $v$. 
By~\eqref{eq:tails}, either $h_i=v$ and $i=k+1$ or $h_i\in T_i$.
If $h_i=v$, then there is nothing to show.
Otherwise, $h_i\in T_i$, so there is a directed path in $D$ from $h_i$ to some $z\in Z$ entirely inside $T_i$.
Since the labels of vertices on  this path decrease by exactly $1$ under $\vp_i$, the path remains directed in $D_i$.
Moreover, $z\in T_i$ and $v\notin T_i$, so $\vp_i(v)-\vp_i(z)=1$.
Thus $\ora{zv}$ is an arc of $D_i$, so $h_i$ can reach $v$.

By~\eqref{eq:tails}, $\dist(t_i)=i-1$, there is a shortest directed path $P$ from $v$ to $t_i$ in $D$.
Every vertex $u$ on $P$ satisfies $\dist(u)<i$, so no vertex on $P$ is in $T_i$.
Therefore, the labels of vertices on $P$ are unchanged under $\vp_i$, so $P$ remains a directed path in $D_i$.

Combining these paths gives a closed directed walk from $t_i$ to $h_i$ to $v$ back to $t_i$ in $D_i$.
This walk is nonempty since $t_i\ne h_i$.
Indeed, $h_i\in T_i$ implies $t_i\notin T_i$ and if $h_i=v$, then $i=k+1$ so $\dist(t_i)=k>0=\dist(v)$.

A shortest nonempty closed directed walk in $D_i$ starting and ending at $t_i$ is a directed cycle through $t_i$.
Let $C_i$ be its underlying  cycle.
As shown above, $C_i$ is good, and $t_i\in V(C_i)\cap(T_{i-1}\setminus T_i)$.

\item \emph{Counting distinct good cycles.}

Suppose $C_i=C=C_j$ for some $i<j$.
Since a good cycle has only one good cyclic traversal, $C_i$ and $C_j$ are traversed the in same direction in $D_i$ and $D_j$.
Thus, for $z\in\{i,j\}$ and for each edge $xy$ on $C_z$ traversed from $x$ to $y$, satisfies $\vp_z(y)-\vp_z(x)=1$ if $x\prec y$ and $\vp_z(y)-\vp_z(x)=1-k$ if $y\prec x$.
Thus $\vp_i(y)-\vp_i(x)=\vp_j(y)-\vp_j(x)$, which implies that $\vp_j-\vp_i$ is constant on the vertices of $C$.

On the other hand, since $i<j$, we have $T_j\subseteq T_i$.
For each vertex $u\in T_i\setminus T_j$ on $C$, we have $\vp_i(u)=\vp(u)-1$ and $\vp_j(u)=\vp(u)$.
Also, $\vp_i(u)=\vp_j(u)$ when $u\in T_j$ or $u\notin T_i$. 
Therefore,
$\vp_j(u)-\vp_i(u)=
\begin{cases}
1,&u\in T_i\setminus T_j,\\
0,&\text{otherwise}.
\end{cases}$
Since $t_i\notin T_i$, $\vp_j(t_i)-\vp_i(t_i)=0$.
However, by~\eqref{eq:tails}, $t_j\in T_{j-1}\setminus T_j\subseteq T_i\setminus T_j$,
so $\vp_j(t_j)-\vp_i(t_j)=1$.
Since $t_i, t_j$ are on $C$, this contradicts that $\vp_j-\vp_i$ is constant on vertices of $C_j$.
\end{enumerate}

Therefore $C_1,\ldots,C_{k+1}$ are
pairwise distinct.
\end{proof}

\section{Counting vertex orderings}
\label{sec:count}

In this section, we bound the number of vertex orderings for which a given cycle is good. 

\begin{lemma}\label{lem:colororders}
For $k\geq 3$, let $H$ be a graph with a vertex $v$ where $\vp:V(H-v)\to\{0,\ldots,k-1\}$ is a proper $k$-coloring.
For each permutation $\pi$ of the $k$ colors, let $\prec_\pi$ be the vertex ordering of $H$ obtained by ordering the vertices of $H-v$ in nondecreasing order of their colors under $\vp$ according to  $\pi$, and placing $v$ last.
Then, each $ak$-cycle in $H$ is good in $H_{\prec_\pi}$ for at most  $2k$ permutations $\pi$.
\end{lemma}
\begin{proof}
Fix an $ak$-cycle $C$, and fix a cyclic traversal $\vec C$. 
We show that at most $k$ permutations $\pi$ make $\vec C$ good. 

$(\ast)$ We first observe that any $k$ consecutive edges of $\vec C$ not involving $v$ must contain a backward edge.
Otherwise, the colors of the $k+1$ vertices on the $k$ edges would be  strictly increasing in the order specified by $\pi$, which contradicts that there are $k$ colors.

\begin{enumerate}[
    wide=\dimexpr\parindent/2\relax,
    label=\emph{{Case \arabic*.}},
]
\item $C$ avoids $v$.

Partition the $ak$ edges of $C$  into $a$ parts, where each part has $k$ consecutive edges on $C$.
By $(\ast)$, each part contains a backward edge, so there are at least $a$ backward edges.
If $\vec C$ is good, then there must be exactly $a$ backward edges.
This further implies that there are  exactly $k-1$ forward edges between  two closest backward edges.

Consequently, the $k$ vertices following each backward edge have all $k$ colors, appearing in the order specified by $\pi$.
In particular, the same ordering of all $k$ colors appears $a$ times around $\vec C$. 

To determine a permutation $\pi$ for which $\vec C$ is good, we only need to choose which of the $k$ colors comes first.
The order of the remaining $k-1$ colors is then uniquely determined by their circular order along $\vec C$.
Thus there are at most $k$ permutations $\pi$ for which $\vec C$ is good.

\item $C$ contains $v$.

Let $\vec C=u_1u_2\ldots u_{ak-1}v$.
Since $v$ is last in $\prec_\pi$, the edge $vu_1$ is backward, while the edge $u_{ak-1}v$ is forward.
If $\vec C$ is good, then at most $a-1$ edges in $E(C)\setminus\{vu_1, u_{ak-1}v\}$ are backward.

Consider the last $(a-1)k$ edges of the path $u_1u_2\ldots u_{ak-1}$.
Partition these edges into $a-1$ parts where each part has $k$ consecutive edges.
By $(\ast)$, each of the $a-1$ parts contains at least one backward edge.
It follows that $u_1u_2, u_2u_3, \ldots,u_{k-2}u_{k-1}$ are forward edges, so $\vp(u_1),\vp(u_2),\ldots,\vp(u_{k-1})$ are distinct and appear in increasing order according to $\pi$. 

These $k-1$ colors determine their relative order in $\pi$. 
There are exactly $k$ possible positions to insert the unique missing color in $\{0, 1, \ldots, k-1\}\setminus\{\vp(u_1), \ldots, \vp(u_{k-1})\}$.
Thus at most $k$ permutations $\pi$ make  $\vec C$ good. 
\end{enumerate}

Since $C$ has two possible cyclic traversals, and $\vec C$ is good for at most $k$ permutations, we conclude $C$ is good for at most $2k$ permutations.
\end{proof}

\section{The conclusion}
\label{sec:conclusion}

We now prove Theorem~\ref{thm:main}.
Given an integer $k\geq 3$, let $G$ be a graph with $\chi(G)>k$. 
Since $\chi(G)>k$, let $H$ be an induced subgraph of $G$ with the minimum number of vertices where $\chi(H)>k$ and every vertex-minimal  subgraph $H'$ of $H$ satisfies $\chi(H')\leq k$. 
Fix $v\in V(H)$.
By minimality, $H-v$ has a proper $k$-coloring $\vp$.
Consider all $k!$ permutations $\pi$ of the $k$ colors and their induced vertex orderings $\prec_\pi$.
For each $\pi$, relabel the colors so that their increasing order agrees with $\pi$.

By Lemma~\ref{lem:orders}, each  $H_{\prec_\pi}$ has at least $k+1$ good cycles.
Since each cycle is good for at most $2k$ permutations by Lemma~\ref{lem:colororders},
$H$ contains at least $\frac{k!(k+1)}{2k}=\frac{(k+1)!}{2k}$
cycles whose lengths are divisible by $k$.
Since every cycle of $H$ is also a cycle of $G$, the desired lower bound follows.

To see that the bound is sharp, consider $K_{k+1}$.
Note that $\chi(K_{k+1})=k+1$.
There are $k+1$ ways to choose the $k$ vertices of a $k$-cycle.
For every choice of $k$ vertices, there are $\frac{(k-1)!}{2}$ distinct $k$-cycles. 
Thus $K_{k+1}$ has exactly
$(k+1)\frac{(k-1)!}{2}
=\frac{(k+1)!}{2k}$
cycles of length divisible by $k$.

\subsection{Further questions}

Our result determines the minimum number of cycles of length $0$ modulo $k$ among graphs with no proper $k$-coloring. 
More generally, Brewster, McGuinness, Moore, and Noel~\cite{BMMN2016} explicitly suggested studying the same question for cycles of length $r$ modulo $k$. 
(See also the discussion in the last section of~\cite{Tuza1992}.)
Let $N_{k,r}(G)$ denote the number of cycles of $G$ with length $r$ modulo $k$.
We have shown that every graph $G$ with no proper $k$-coloring satisfies $N_{k, 0}(G)\geq N_{k, 0}(K_{k+1})$.
As mentioned in the introduction,  Moore and West~\cite{MooreWest2021} proved the analogous inequality for $r=1$.
For $r=2$, the inequality holds since $N_{k, 2}(K_{k+1})=0$. 
Thus, it is natural to ask if this phenomenon continues: 

\begin{question}
Let $k\ge4$ and $r\in\{3, \ldots,  k-1\}$.
Does every graph $G$ with no proper $k$-coloring satisfy
$N_{k,r}(G)\ge
N_{k,r}(K_{k+1})$?
\end{question}

Another question is to characterize equality in Theorem~\ref{thm:main}.
It is natural to restrict attention to critical graphs.
Recall that a graph $G$ is \emph{$k$-critical} if $\chi(G)=k$ and every proper subgraph of $G$ has a proper  $(k-1)$-coloring. 
Kim and Picollelli~\cite{KimPicollelli2024}
showed that $K_4$ is the only $4$-critical graph with exactly four cycles of length $0$ modulo $3$.
We ask whether this extends to larger values of $k$. 

\begin{question}
Let $k\ge4$.
Is $K_{k+1}$ the unique $(k+1)$-critical graph with exactly $\frac{(k+1)!}{2k}$ cycles of length $0$ modulo $k$?
\end{question}

\bigskip
\noindent\textbf{Declaration of LLM use.}
\par\smallskip
\noindent
The authors used OpenAI's models to assist with mathematical exploration and argument refinement. 
The authors developed and verified the mathematical arguments, revised the exposition, and take full responsibility for the contents of the article.


\end{document}